\documentclass[reqno, oneside]{amsart}
\usepackage{style}

\begin{document}
\title[Toeplitz operators]{Schatten properties of Toeplitz operators \\ on the Paley-Wiener space}
\author{R.~V.~Bessonov}

\address{St.Petersburg State University ({\normalfont \hbox{29b}, 14th Line V.O., 199178, St.Petersburg,  Russia}) and St.Petersburg Department of Steklov Mathematical Institute of Russian Academy of Science ({\normalfont 27, Fon\-tan\-ka, 191023, St.Petersburg, Russia})}
\email{bessonov@pdmi.ras.ru}

\thanks{The work is supported by RFBR grant mol\_a\_dk 16-31-60053, by Grant MD-5758.2015.1, and by ``Native towns'', a social investment program of PJSC ``Gazprom Neft'' }

\begin{abstract}
We collect several old and new descriptions of Schatten class Toeplitz operators 
on the Paley-Wiener space and answer a question on discrete Hilbert transform commutators posed by Richard Rochberg.
\end{abstract}

\maketitle

\section{Introduction}\label{s1}
Given a bounded function $\phi$ on the real line, $\R$, consider the Toeplitz operator~$T_\phi$ on the classical Paley-Wiener space $\pw_{a}$, 
\begin{equation}\label{eq1}
T_\phi\colon f \mapsto P_{a} (\phi f), \qquad f \in \pw_{a}.
\end{equation}
The space $\pw_{a}$ could be regarded as the subspace in $L^2(\R)$ of functions with Fourier spectrum in the interval $[-a,a]$, symbol $P_{a}$ above denotes the orthogonal projection in $L^2(\R)$ to $\pw_{a}$. Basic theory of Toeplitz operators on $\pw_a$ can be found in paper \cite{Roch87} by R.\ Rochberg. 

\medskip

We are interested in description of Schatten class Toeplitz operators on 
$\pw_{a}$ in terms of their standard symbols. By the standard symbol of an operator in \eqref{eq1} we mean the entire function 
$\phi_{st} = \F^{-1}\chi_{2a}\,\F\phi$, where $\F$ denotes the Fourier transform on the Schwartz space of tempered distributions, and $\chi_{2a}$ is the indicator function of the interval $(-2a,2a)$.  As we will see, a Toeplitz operator $T_\phi$ on $\pw_a$ belongs to the Schatten class $\Ss^p$, $0<p<\infty$, if and only if $e^{2iax}\phi_{st}$ belongs to a discrete oscillation Besov space introduced in 1987 by R.\ Rochberg \cite{Roch87}. Its definition we now recall.

\medskip
For a measure $\mu$ on $\R$ and a function $f \in L^{1}_{{\rm loc}}(\mu)$, the oscillation of order $n$ of $f$ on an interval $I \subset \R$ with respect to $\mu$  is defined by
$$
\osc(f, I, \mu, n) =  \inf_{P_n}\frac{1}{\mu(I)}\int_{I}|f(x) - P_n(x)|\,d\mu(x), 
$$ 
where the infimum is taken over all polynomials $P_n$ of degree at most $n$. If $\mu(I) = 0$, we put $\osc_{I}(f, I, \mu, n)  = 0$. Define the family $\I_a$ of closed intervals 
$$
I_{a,j,k} = \left[\frac{2\pi}{a}\, k \,2^j, \frac{2\pi}{a}(k+1)2^j\right], \qquad j,k \in \Z, \quad j\ge0.
$$  
Note that endpoints of intervals in $\I_a$ belong to the lattice  $\Z_{a} = \bigl\{\frac{2\pi}{a}k, \; k \in \Z\bigr\}$. Let $p$ be a positive real number, and let $[\tfrac{1}{p}]$ be the integer part of $\tfrac{1}{p}$. The discrete oscillation Besov space $\Bb_p(a, \osc) = \Bb_{p,p}^{1/p}(\Z_a,\mu_a,\osc)$ is defined by
$$
\Bb_p(a,\osc) = \left\{f \in L^{1}_{\rm{loc}}(\mu_a): \;\; \|f\|_{\Bb_p(a,\osc)}^{p} = \sum\nolimits_{I \in \I_a}\osc\left(f, I, \mu_a, \bigl[\tfrac{1}{p}\bigr]\right)^{p} <\infty \right\},
$$ 
where $\mu_a = \frac{2\pi}{a}\sum_{x \in \Z_a} \delta_{x}$ is the normalized counting measure on $\Z_a$. 

\medskip

Our main result is the following theorem.
\begin{Thm}\label{t1}
Let $a$, $p$ be positive real numbers, let $\phi$ be a bounded function on $\R$, and let $\phi_{st}$ be the standard symbol of the Toeplitz operator $T_\phi$ on $\pw_a$. Then we have $T_\phi\in\Ss^p(\pw_a)$ if and only if $e^{2iax}\phi_{st} \in \Bb_p(4a,\osc)$. Moreover, $\|T_\phi\|_{\Ss^p}$ is comparable to $\|e^{2iax}\phi_{st}\|_{\Bb_p(4a,\osc)}$ with  constants depending only on $p$.
\end{Thm}
Theorem \ref{t1} complements a classical description of Toeplitz operators in $\Ss^p(\pw_a)$ given by R.~Rochberg \cite{Roch87} for $1 \le p<\infty$ and extended by V.~Peller \cite{Peller88} to the whole range $0<p<\infty$.  To formulate the result, consider a system $\{\nu_{j}\}_{j \le -1}$ of infinitely smooth functions on $\R$ such that 
$$
\supp \nu_j \subset [2^{j-1}, 2^{j}], \quad 0 \le \nu_j \le 1,  \quad \nu_{j-1}(x) = \nu_{j}(x/2), \quad \sum \nu_j = 1 \mbox{ on }(0, \tfrac{1}{3}].
$$ 
Define $\nu_{j}(x) = \nu_{-j}(1 - x)$ for real $x \ge \frac{1}{2}$ and integer $j \ge 1$, put $\nu_0 = 1 - \sum_{j \neq 0} \nu_j$ for $j=0$. Finally, let $\nu_{a,j}(x) = \nu_j((x+a)/2a)$ for all $x \in [-a,a]$ and $j \in \Z$. Observe that system $\{\nu_{a,j}\}$ provides a resolution of unity on the interval $[-a,a]$ by functions supported on subintervals $I_j$ whose lengths are comparable to the distance from $I_j$ to the endpoints of $[-a,a]$. Rochberg-Peller theorem says that $T_\phi$ is in $\Ss^p(\pw_a)$ for $0<p<\infty$ if and only if
$$
a\sum_{j \in \Z} 2^{-|j|}\cdot \|\F^{-1}(\nu_{2a,j}\cdot \F \phi)\|_{L^p(\R)}^{p} < \infty,
$$
with control of the norms. R.~Rochberg gives yet another characterization of Toeplitz operators in class $\Ss^p(\pw_a)$, $1 \le p<\infty$, in terms of a reproducing kernel decomposition of their standard symbols, see Theorem~5.3 in \cite{Roch87}. Both the statement and the proof of his result for $p=1$ contain errors that we correct in Section~\ref{s2}. 

\medskip

As a consequence of Theorem~\ref{t1}, we obtain the following result. 
\begin{Thm}\label{t2}
Let $a>0$. The discrete Hilbert transform commutator
$$
C_\psi: f \mapsto \frac{1}{\pi}\dashint_{\Z_{a}}\frac{\psi(x) - \psi(t)}{x -t}f(t)\,d\mu_a(t), \qquad f \in L^2(\mu_a),
$$ 
belongs to the trace class $\Ss^1(L^2(\mu_a))$ if and only if $\psi \in \Bb_1(a,\osc) \cap L^\infty(\Z_{a})$. 
\end{Thm}
This answers the question posed by R.~Rochberg in 1987. See Section~\ref{s4} for a summary of results on discrete Hilbert transform commutators and an analogue of Theorem \ref{t2} for the case $0<p<1$. 

\medskip

We would like to mention papers \cite{Tor91}, \cite{Tor98} by R.~Torres for readers interested in wavelet characterizations and interpolation theory of discrete Besov spaces. The problem of membership in Schatten classes $\Ss^p$ for general truncated Toeplitz operators has been recently studied by P.~Lopatto and R.~Rochberg \cite{RochLop16}, see also Section~4.3 in author's paper \cite{Bes15}. 

\medskip

\section{Proof of Theorem \ref{t1} for $1<p<\infty$}
Theorem \ref{t1} for $1<p<\infty$ follows from known results. Let $\Bb_p(\R) = \dot{\Bb}_{p,p}^{1/p}(\R)$ be the standard homogeneous Besov space on the real line~$\R$, see, e.g., Chapter 3 in \cite{Peetre76} for definition and basic properties. Given a Toeplitz operator $T_{\phi}$ on $\pw_a$ with symbol $\phi\in L^\infty(\R)$, we denote 
$$
\phi_{st}^{-} = \F^{-1}\chi_{(-2a,0)}\F\phi, \qquad \phi_{st}^{+} = \F^{-1}\chi_{[0,2a)}\F\phi,
$$ 
where $\chi_{S}$ is the indicator function of a set $S$. As usual, $\F$ stands for the Fourier transform on the Schwartz space of tempered distributions. The following result is a combination of Theorem 5.1  and its Corollary in \cite{Roch87}.
\begin{NMT}[R.~Rochberg]
Let $1 < p < \infty$ and let $a>0$. Then a Toeplitz operator $T_{\phi}$ on $\pw_a$ belongs to $\Ss_p(\pw_a)$ if and only if 
$\|e^{2iax}\phi_{st}^{-}\|_{\Bb_p(\R)} + \|e^{-2iax}\phi_{st}^{+}\|_{\Bb_p(\R)}$ is finite, in which case $\|T_{\phi}\|_{\Ss^p}$ is comparable to 
$\|e^{2iax}\phi_{st}^{-}\|_{\Bb_p(\R)} + \|e^{-2iax}\phi_{st}^{+}\|_{\Bb_p(\R)}$ with constants depending only on $p$.
\end{NMT}
Denote by $\E_{a}$ the set of tempered distributions whose Fourier transforms are supported on the interval $[-a,a]$. Next result is Theorem 1 in \cite{Tor98}. 
\begin{NMT}[R.~Torres]
Let $1 < p < \infty$ and let $f$ be a function in $\E_{a} \cap \Bb_{p}(\R)$ for some $a>0$. Then its restriction to $\Z_{2a}$ belongs to $\Bb_{p}(2a,\osc)$ and $\|f\|_{\Bb_{p}(2a,\osc)}$ is comparable to $\|f\|_{\Bb_{p}(\R)}$ with constants depending only on $p$.
Moreover, every sequence in $\Bb_p(a,\osc)$ is the restriction to $\Z_a$ of a unique function (modulo polynomials) in $\E_a\cap\Bb_p(\R)$.   
\end{NMT}

\noindent{\bf Proof of Theorem \ref{t1} ($1<p<\infty$).} Let $\phi$ be a bounded function of $\R$ and let 
$\phi_{st} = \F^{-1}\chi_{(-2a,2a)}\F\phi$ be the standard symbol of the Toeplitz operator  $T_{\phi} \in \Ss^p(\pw_a)$. Then functions
$e^{2iax}\phi_{st}^{-}$, $e^{-2iax}\phi_{st}^{+}$ belong to $\E_{2a} \cap \Bb_{p}(\R)$ by R.~Rochberg's theorem above. From theorem by 
R.~Torres we see that $e^{2iax}\phi_{st}^{-} \in \Bb_p(4a,\osc)$ and $e^{-2iax}\phi_{st}^{+} \in \Bb_p(4a,\osc)$ with control of the norms. Now observe that $e^{4iax} = 1$  and $e^{2iax}\phi_{st} = e^{2iax}\phi_{st}^{-} + e^{-2iax}\phi_{st}^{+}$ on $\Z_{4a}$, hence $e^{2iax}\phi_{st} \in \Bb_{p}(4a,\osc)$.   

\medskip

Conversely, assume that the restriction of $e^{2iax}\phi_{st}$ to $\Z_{4a}$ is in $\Bb_{p}(4a,\osc)$. Using theorem by R.~Torres, find a function $f \in \E_{2a}\cap\Bb_p(\R)$ such that its restriction to $\Z_{4a}$ agrees with $e^{2iax}\phi_{st}$. Put $f^{-} = \F^{-1}\chi_{(-2a,0)}\F f$ and $f^{+} = \F^{-1}\chi_{[0,2a)}\F f$. Observe that 
$\tilde \phi = e^{-2iax}f^{+} + e^{2iax}f^{-}$ is an entire function of exponential type at most $2a$ coinciding with $\phi_{st}$ on $\Z_{4a}$. 
Since $\phi_{st}$, $\tilde\phi$ are the first order distributions supported on the finite interval $[-2a,2a]$, we have $|\tilde\phi(x)| + |\phi(x)| \le c+ c|x|$ for all $x\in \R$ and a constant $c\ge0$. It follows that the entire function
$\frac{\tilde\phi - \phi}{z}$ of exponential type at most $2a$ is bounded on $\R$ and vanishes on $\Z_{4a}\setminus\{0\}$, hence $\tilde\phi - \phi_{st} = p\sin(2az)$ for a polynomial $p$ of degree at most $1$. Therefore, we have $T_{\phi} = T_{\phi_{st}} = T_{\tilde\phi}$ on  $\pw_a$, see Section~2.D in \cite{Roch87}.
Since $f^{\pm}\in \Bb_p(\R)$, we can use R.~Rochberg's theorem  and conclude 
that $T_{\tilde\phi} \in \Ss^p(\pw_a)$ with control of the norms: $\|T_{\tilde\phi}\|_{\Ss^p}$ is controllable by $\|e^{2iax}\tilde\phi^{-}\|_{\Bb_p(\R)} + \|e^{-2iax}\tilde\phi^{+}\|_{\Bb_p(\R)}  \le c_p \|f\|_{\Bb_p(\R)} \le \tilde c_p \|e^{2iax}\phi_{st}\|_{\Bb_p(4a,\osc)}$. 
 \qed

\medskip

\section{Reproducing kernel decomposition of standard symbols}\label{s2}
In this section we show that the standard symbol of a Toeplitz operator on $\pw_a$ from class $\Ss^p$ could be represented as a linear combination of normalized reproducing kernels of $\pw_{2a}$ with coefficients $c_k$ such that $\sum |c_k|^p < \infty$. We consider only the case $0<p\le1$. Proposition \ref{p1} below is a corrected version of Theorem 5.3 in \cite{Roch87}. In the original statement the author of \cite{Roch87} forgot to normalize the exponentials in formula $(5.6)$ of \cite{Roch87}. More importantly, he used the fact that the Fourier multiplier $f \mapsto \F^{-1}\chi_{[0,1]}\F f$ is bounded on $\Bb_{p}(\R)$. This is not the case for $p=1$. Here is a more accurate implementation of the ideas from \cite{Roch87}.   

\medskip

Let $\psi$ be a bounded function on the real line $\R$. Consider the standard Hardy space $H^2$ in the upper half-plane $\C^+ = \{\lambda \in \C: \;\; \Im\lambda >0\}$ of the complex plane~$\C$. Denote by $H^2_-$ the anti-analytic subspace $\{f \in L^2(\R): \; \bar f \in H^2\}$ of $L^2(\R)$. Recall that the classical Hankel operator $H_{\psi}: H^2 \to H^2_-$ is defined by
$$
H_{\psi}: f \mapsto P_{-}(\psi f), \qquad f \in H^2,
$$
where $P_-$ denotes the orthogonal projection from $L^2(\R)$ to $H^2_-$. The operator $H_{\psi}$ is completely determined by its standard anti-analytic symbol $\psi_{st} = \F^{-1}\chi_{(-\infty, 0)}\F\psi$. The latter means that $H_{\psi}f = H_{\psi_{st}}f$ for all  $f \in H^2$ such that $\sup_{x \in \R} |xf(x)| < \infty$. Take a positive number $\eps>0$ and define the sets $\U_{\eps}^{+}$, $\U_{\eps}^{-}$ by
$$
\U_{\eps}^{\pm} = \{\lambda \in \C: \;\; \lambda = (1+\eps)^m(\eps x \pm i); \;\; x,m \in \Z\}.
$$
For $\lambda \in \C^+$ let us denote by $k_{\lambda}$ the reproducing kernel of $H^2$ at $\lambda$, $k_\lambda = -\frac{1}{2\pi i}\frac{1}{z - \bar\lambda}$. 
\begin{NMT}[R.~Rochberg \cite{Roch85}]
There exists a number $\eps>0$ such that $H_{\psi} \in \Ss^p(H^2)$ if and only if $\psi_{st} = \sum_{\lambda \in \U_{\eps}^{+}} c_\lambda \frac{\ov{k_{\lambda}}}{\|k_{\lambda}\|^2}$, where $\sum|c_\lambda|^p$ is finite and the infimum of $\sum|c_\lambda|^p$ over all possible representations of $\psi_{st}$ in this form is comparable to $\|H_{\psi}\|_{\Ss^p}^{p}$ with constants depending only on $p\in(0,\infty)$.
\end{NMT}
Remark that for $p\in (0,1]$ the series defining $\psi_{st}$ in the theorem above converges absolutely to a bounded function on $\R$, while for $p>1$ the convergence holds only in the Besov space $\Bb_p(\R)$ (one need to extract constant terms from every summand to get the convergent series, see discussion in \cite{Roch85}). 
In order to prove an analogous result for Toeplitz operators on the Paley-Wiener space, let us consider the sets  
\begin{align*}
\U_{\eta a, \eps}^{\pm} = \{\lambda \in \U^{\pm}_{\eps}: |\Im\lambda| > \tfrac{\eps}{\eta a}\}, \qquad 
\Lambda_{\eta a,\eps} = \U_{\eta a,\eps}^{-} \cup \Z_{\eta a} \cup \U_{\eta a,\eps}^{+}.
\end{align*}
Here $\Z_{\eta a} = \{\frac{2\pi}{\eta a}k, \; k\in \Z\}$. Next, for $a>0$ and $\lambda \in \C$, denote by $\rho_{a,\lambda}$ the reproducing kernel of the space $\pw_a$ at the point $\lambda$. Recall that
$$
\rho_{a,\lambda}: z \mapsto \frac{1}{\pi}\frac{\sin a(z-\bar\lambda)}{z-\bar\lambda}, \qquad z \in \C.
$$ 
We are going to prove the following proposition.
\begin{Prop}\label{p1}
Let $a>0$ and let $\phi \in L^\infty(\R)$. There exist $\eps>0$, $\eta>1$ such that $T_{\phi}\in\Ss^p(\pw_a)$ if and only if 
$\phi_{st} = \sum_{\lambda \in \Lambda_{\eta a,\eps}} c_\lambda \frac{\rho_{2a,\lambda}}{\|\rho_{a,\lambda}\|^2}$, where $\sum_{\lambda}|c_\lambda|^{p}$ is finite and the infimum of $\sum|c_\lambda|^p$ over all possible representations of $\phi_{st}$ in this form is comparable to
$\|T_{\phi}\|_{\Ss^p}^{p}$ with constants depending only on $p \in (0,1]$.
\end{Prop}
We will show how to reduce Proposition \ref{p1} to the above theorem for Hankel operators using a splitting of the standard symbol into three pieces: analytic, anti-analytic and a piece with ``small'' Fourier support.  

\medskip

The following two results for $0<p\le1$ are consequences of Lemma~1  and Lemma~2 from \cite{Peller88}. The range $1 \le p < \infty$ has been treated earlier in \cite{Roch87}, see also Section~2 in \cite{Smith}. 
\begin{Lem}\label{l1}
Let $a>0$ and let $\phi \in L^\infty(\R)$. There exist bounded functions $\phi_\ell$, $\phi_c$, and $\phi_{r}$ such that  $T_{\phi} = T_{\phi_{\ell}} + T_{\phi_{c}} + T_{\phi_{r}}$ on $\pw_a$,
$$
\supp\F\phi_\ell \subset [-4a, -\tfrac{a}{2}], \quad \supp\F\phi_c \subset [-a, a], \quad \supp\F\phi_r \subset [\tfrac{a}{2}, 4a],
$$
and we have $\|T_{\phi_{s}}\|_{\Ss^p} \le c_p \|T_{\phi}\|_{\Ss^p}$ for every $s = \ell, c, r$ for $T_{\phi} \in \Ss^p(\pw_a)$. Here $c_p$ is a constant depending only on $p$.
\end{Lem}
\begin{Lem}\label{l2}
Let $a>0$ and let $\phi \in L^\infty(\R)$ be such that $\supp\hat\phi\subset[-a,a]$. Then we have $T_\phi\in\Ss^p(\pw_a)$ if and only if $\phi \in L^p(\R)$, in which case $\|\phi\|_{L^p(\R)}$ is comparable to $\|T_\phi\|_{\Ss^p}$ with constants depending only on $p$.
\end{Lem}

\medskip

\noindent{\bf Proof of Proposition \ref{p1}.} Let $\phi\in L^\infty(\R)$ and let $\phi_{st} =  \F^{-1}\chi_{(-2a,2a)}\,\F\phi$ be the standard symbol of the operator $T_{\phi}$ on $\pw_a$. Then $T_\phi = T_{\phi_{st}}$, see Section~2.D in \cite{Roch87}.  Suppose that $\phi_{st} = \sum_{\lambda \in \Lambda_{\eta a,\eps}} c_\lambda \frac{\rho_{2a,\lambda}}{\|\rho_{a,\lambda}\|^2}$ for some $\eps>0$, $\eta>0$, and some coefficients $c_\lambda$ such that $\sum_{\lambda\in \Lambda_{\eta a,\eps}}|c_\lambda|^{p}<\infty$. It follows from the estimate 
$$
\frac{|\rho_{2a,\lambda}(z)|}{\|\rho_{a,\lambda}\|^2} \le c e^{2a|\Im z|}, \qquad z \in \C, \quad \lambda \in \C,
$$  
that this series converges absolutely to an entire function of exponential type at most $2a$ bounded on the real line $\R$.
By triangle inequality (see, e.g.,\ Theorem A1.1 in \cite{PeBook}), we have
$$
\|T_{\phi}\|_{\Ss^p}^{p} =\|T_{\phi_{st}}\|_{\Ss^p}^{p} \le \Bigl(\sum_{\lambda\in \Lambda_{\eta a,\eps}}|c_\lambda|^{p}\Bigr)\sup_{\lambda \in \C}\|T_{\phi_{\lambda}}\|_{\Ss^p}^{p},
$$
where we denoted $\phi_{\lambda} = \frac{\rho_{2a,\lambda}}{\|\rho_{a,\lambda}\|^2}$. Take $\lambda\in \C$. For every $f,g \in \pw_a$ we have 
$$
(T_{\rho_{2a,\lambda}} f,g) = (f \bar g, \rho_{2a,\bar \lambda}) = f(\bar\lambda)\cdot \ov{g(\lambda)} = (f, \rho_{a,\bar\lambda}) (\rho_{a, \lambda},g). 
$$
It follows that the operator $T_{\phi_{\lambda}}$ has rank one and $\|T_{\phi_\lambda}\|_{\Ss^p} = 1$. Hence $T_{\phi}$ belongs to $\Ss^p(\pw_a)$ and
$\|T_{\phi}\|_{\Ss^p}^{p} \le \sum_{\lambda}|c_\lambda|^p$.

\medskip

Now let $\phi$ be a bounded function on $\R$ such that $T_{\phi}\in\Ss^p(\pw_a)$. We want to show that the standard symbol $\phi_{st} =  \F^{-1}\chi_{(-2a,2a)}\,\F\phi$ of $T_\phi$ can be represented in the form 
$$
\phi_{st} = \sum_{\lambda \in \Lambda_{\eta a,\eps}} c_\lambda \frac{\rho_{2a,\lambda}}{\|\rho_{a,\lambda}\|^2}
$$ 
for some positive numbers $\eps$, $\eta$ depending only on $p$ and a sequence $\{c_\lambda\}$ such that $\sum_{\lambda}|c_\lambda|^{p}$ is comparable to $\|T_{\phi}\|_{\Ss^p}^{p}$. By Lemma \ref{l1}, it suffices to consider separately the following three cases: $(1)$ $\supp\hat\phi \subset (-\infty, 0]$; $(2)$ $\supp\hat\phi \subset [-a, a]$; $(3)$ $\supp\hat\phi \subset [0, +\infty)$. Let us treat the third case first. Denote by $M_{e^{-iax}}$ the operator of multiplication by $e^{-iax}$ on~$L^2(\R)$. Since $\supp\hat\phi \subset [0,+\infty)$, we have 
$$
H_{e^{-2iax}\phi} = M_{e^{-iax}}T_{\phi}P_{a}M_{e^{-iax}},
$$
where $H_{e^{-2iax}\phi}: H^2 \to H^{2}_{-}$ is the Hankel operator with symbol $\psi = e^{-2iax}\phi$. In particular,  we have $\|H_\psi\|_{\Ss^p} \le \|T_{\phi}\|_{\Ss_p}$. By Rochberg's Theorem above, the anti-analytic function $\psi_{st} = \F^{-1}\chi_{(-\infty,0)}\F e^{-2iax}\phi$ admits the following representation:
$$
\psi_{st} = \sum_{\lambda \in \U_{\eps}^{+}} c_\lambda \frac{\ov{k_{\lambda}}}{\|k_{\lambda}\|^2}, 
$$
where $\sum_{\lambda \in \U_{\eps}^{+}}|c_{\lambda}|^p$ is comparable to $\|H_{\psi}\|^{p}_{\Ss^p}$, and $\eps>0$ does not depend on $\psi$. This gives us decomposition for $\phi_{st}$: 
$$
\phi_{st} = e^{2iax}\psi_{st} = \sum_{\lambda \in \U_{\eps}^{+}} c_\lambda \frac{e^{2iax}\ov{k_{\lambda}}}{\|k_{\lambda}\|^2} 
= \sum_{\lambda \in \U_{\eps}^{+}} c_\lambda \frac{P_{2a}(e^{2iax}\ov{k_{\lambda}})}{\|k_{\lambda}\|^2}, 
$$
where $P_{2a}$ denotes the orthogonal projection in $L^2(\R)$ to $\pw_{2a}$. It is easy to see that 
$P_{2a}(e^{2iax}\ov{k_{\lambda}}) = e^{2ia\bar\lambda}\rho_{2a,\bar\lambda}$ and $\|\rho_{a,\bar\lambda}\|^2 \le 2e^{2a\Im\lambda}\cdot\|k_{\lambda}\|_{L^2(\R)}^{2}$, hence
\begin{equation}\notag
\phi_{st} = \sum_{\lambda \in \U_{\eps}^{-}} c_{\bar\lambda} \beta_\lambda \frac{\rho_{2a,\lambda}}{\|\rho_{a,\lambda}\|^2}
\end{equation}
for some complex numbers $\beta_\lambda$ such that $\sup_{\lambda}|\beta_\lambda| \le 2$. Next, in the case where $\supp\phi \subset (-\infty,0]$ we can consider the adjoint operator $T^\ast_{\phi} = T_{\phi_{st}^{*}}$ with the standard symbol $\phi_{st}^{*}: z \mapsto \ov{\phi_{st}(\bar z)}$ and conclude that in this situation
\begin{equation}\notag
\phi_{st} = \sum_{\lambda \in \U_{\eps}^{+}} \ov{c_\lambda\beta_{\bar\lambda}} \frac{\rho_{2a,\lambda}}{\|\rho_{a,\lambda}\|^2}.
\end{equation}
Now let $\supp\phi \subset [-a,a]$. By Lemma \ref{l2}, we have $\phi \in L^p(\R)$. In particular, $\phi \in \pw_{2a}$ and 
Plancherel-Polya theorem \cite{PP} yields the following decomposition:
$$
\phi = \phi_{st} = \frac{\pi}{2a} \sum_{\lambda \in \Z_{2a}} f(\lambda)\rho_{2a,\lambda}, \qquad \sum_{\lambda \in \Z_{2a}}|f(\lambda)|^p \le c_p a^p \|\phi\|_{L^p(\R)}^{p},
$$
where the constant $c_p$ depends only on $p$. Put $\Lambda_{\eps} = \U^{+}_{\eps} \cup \Z_{2a} \cup \U^{-}_{\eps}$. To summarize, we have proved that for every bounded function $\phi$ on $\R$ such that $T_{\phi} \in \Ss^p(\pw_{a})$ there are coefficients $c_{\lambda}$, $\lambda \in \Lambda_{\eps}$, such that
\begin{equation}\label{eq3}
\phi_{st} = \sum_{\lambda \in \Lambda_{\eps}}c_{\lambda} \frac{\rho_{2a,\lambda}}{\|\rho_{a,\lambda}\|^2}, 
\qquad \sum_{\lambda \in \Lambda_{\eps}}|c_{\lambda}|^p \le c_p \|T_{\phi}\|_{\Ss^p}^{p}.
\end{equation}
It remains to show that the set $\Lambda_{\eps}$ and coefficients $c_{\lambda}$ in this decomposition could be be replaced by the set $\Lambda_{\eta a,\eps}$  and some new coefficients $c_{\lambda}$ satisfying the second estimate in \eqref{eq3}. To this end, for every point $\lambda \in \Lambda_{\eps}$ denote by $\zeta_\lambda$ the nearest point to $\lambda$ in $\Lambda_{\eta a,\eps} \subset \Lambda_{\eps}$, where $\eta = 2^{k}$ and  $k\in \Z$ is a positive integer number that will be specified later. Consider the function 
$$
\tilde \phi^{(1)} = \sum_{\lambda \in \Lambda_{\eps}}c_{\lambda} \frac{\rho_{2a,\zeta_\lambda}}{\|\rho_{a,\zeta_\lambda}\|^2}
=\sum_{\lambda \in \Lambda_{\eta a,\eps}}\tilde c_{\lambda}^{(1)} \frac{\rho_{2a,\lambda}}{\|\rho_{a,\lambda}\|^2}, 
\qquad \tilde c_{\lambda}^{(1)} = \sum_{\nu \in \Lambda_{\eps},\; \zeta_\nu = \lambda}c_\nu
$$ 
Note that $\tilde\phi^{(1)}$ has the required representation and $\sum|\tilde c_{\lambda}^{(1)}|^p \le \sum|c_{\lambda}|^p$. Moreover, we have
$
\|T_{\phi} - T_{\tilde\phi^{(1)}}\|^{p}_{\Ss^p} \le \sum_{\lambda \in \Lambda_{\eps} \setminus \Lambda_{\eta a,\eps}}|c_\lambda|^p \cdot \|T_{\phi_\lambda} - T_{\phi_{\zeta_\lambda}}\|_{\Ss^p}^{p}.
$ 
On the other hand, the quasi-norm in $\Ss_{p}$ of the rank two operator 
$$
T_{\phi_{\lambda}} - T_{\phi_{\zeta_\lambda}} = \frac{\rho_{a,\lambda}}{\|\rho_{a,\lambda}\|} \otimes \frac{\rho_{a,\lambda}}{\|\rho_{a,\bar\lambda}\|}-
\frac{\rho_{a,\zeta_\lambda}}{\|\rho_{a,\zeta_\lambda}\|} \otimes \frac{\rho_{a,\zeta_\lambda}}{\|\rho_{a,\zeta_\lambda}\|}
$$
does not exceed 
$$
2^{\frac{1}{p}}\left\|\frac{\rho_{a,\zeta_\lambda}}{\|\rho_{a,\zeta_\lambda}\|} - \frac{\rho_{a,\lambda}}{\|\rho_{a,\lambda}\|}\right\|_{L^2(\R)} 
\le 2^{\frac{1}{p}+\frac{1}{2}}\left(1 - \frac{\Re\rho_{a,\zeta_\lambda}(\lambda)}{\|\rho_{a,\zeta_\lambda}\|\cdot\|\rho_{a,\lambda}\|}\right)^{\frac{1}{2}}.
$$ 
Since $|\zeta_\lambda - \lambda| \le \frac{2\pi}{\eta a}$ for all $\lambda$ by construction, one can choose a large number $\eta = 2^k$ so that 
$\|T_{\phi} - T_{\tilde\phi}\|_{\Ss^p}^{p} \le \frac{1}{2}\|T_{\phi}\|_{\Ss^p}^{p}$. Clearly, this choice of $\eta$ does not depend on $\phi$ and $a$. 
Iterating the process, we see that there are functions 
$$
\tilde\phi^{(n)}=\sum_{\lambda \in \Lambda_{\eta a,\eps}}\tilde c_{\lambda}^{(n)} \frac{\rho_{2a,\lambda}}{\|\rho_{a,\lambda}\|^2}, \quad n= 1,2, \ldots,
$$ 
such that $\|T_{\phi} - T_{\tilde\phi^{(1)}} - \ldots T_{\tilde\phi^{(n)}}\|_{\Ss^p}^{p} \le \frac{1}{2^{n}}\|T_{\phi}\|_{\Ss^{p}}^{p}$ and
$\sum_{n,\lambda} |\tilde c_{\lambda}^{(n)}|^{p} \le c_p^p \|T_{\phi}\|_{\Ss^{p}}^{p}$. Since $\Ss^{p}(\pw_a)$ is a complete quasi-normed space and a Toeplitz operator on $\pw_a$ is zero if and only if its standard symbol is zero (see Section 2.D in \cite{Roch87}), this gives us the required decomposition of $\phi_{st}$ with coefficients $c_\lambda = \sum_{n\ge 1} \tilde c_{\lambda}^{(n)}$, $\lambda \in \Lambda_{\eta a,\eps}$.

\medskip


\section{Interpolation of discrete Besov sequences}\label{s3}
Denote by $\pw_{[0,a]}$ the Paley-Wiener space of functions in $L^2(\R)$ with Fourier spectrum in the interval $[0,a]$. Recall that the reproducing kernel 
$k_{a,\lambda}$ of the space $\pw_{[0,a]}$ at a point $\lambda \in \C_+$ has the form
$$
k_{a,\lambda}(z) = -\frac{1}{2\pi i}\frac{1-e^{ia(z-\bar\lambda)}}{z-\bar\lambda}, \qquad z \in \C.
$$
Denote by $\Cc_0(\Z_a)$ the set of functions on $\Z_a$ tending to zero at infinity. 
Our aim in this section is to prove the following proposition.
\begin{Prop}\label{P4}
Let $0<p\le 1$, let $\Lambda$ be the set $\Lambda_{\eta a, \eps}$ from Proposition \ref{p1}, and let $F = \sum_{\lambda \in \Lambda} c_\lambda \frac{k_{a,\lambda}}{\|k_{\frac{a}{2},\lambda}\|^2}$ for some $c_\lambda \in \C$ such that $\sum_{\lambda \in \Lambda}|c_{\lambda}|^{p} < \infty$. Then the restriction of $F$ to $\Z_a$ belongs to $\Bb_{p}(a, \osc) \cap \Cc_0(\Z_{a})$. Conversely, for every function $f \in \Bb_{p}(a, \osc)$ there exists the unique function $F$ as above and a polynomial $q$ of degree at most $[\frac{1}{p}]$ such that $f = q + F$ on $\Z_a$. Moreover, the infinum of $\sum_{\lambda \in \Lambda}|c_{\lambda}|^{p}$ over all possible representations of $F=\sum_{\lambda \in \Lambda} c_\lambda \frac{k_{a,\lambda}}{\|k_{\frac{a}{2},\lambda}\|^2}$ in this form is comparable to $\|f\|_{\Bb_p(\osc,a)}^{p}$ with constants depending only on $p$.  
\end{Prop}
The proof of Proposition \ref{P4} is based on the following lemma.
\begin{Lem}\label{l17}
We have $\|k_{a,\lambda}\|_{\Bb_p(a,\osc)} \le c_p \|k_{\frac{a}{2},\lambda}\|^{2}$ for every $a>0$, $0<p\le 1$, and $\lambda \in \C$, where the constant $c_p$ depends only on $p$.
\end{Lem}
\beginpf At first, consider the points $\lambda$ in the support of $\mu_a$. For $\lambda \in \Z_a$ we have
$$
k_{a,\lambda}(x) = 
\begin{cases}
\|k_{a,\lambda}\|^2, \quad &x = \lambda;\\
0, &x \in \supp\mu_a \setminus\{\lambda\}.
\end{cases}
$$ 
Taking $P_I = 0$ for all intervals $I \in \I_a$ in the definition of $\osc\bigl(k_{a,\lambda}, I, \mu_a, \bigl[\tfrac{1}{p}\bigr]\bigr)$, we obtain the estimate
\begin{align*}
\|k_{a,\lambda}\|_{\Bb_p(a,\osc)}^{p} 
&\le \sum_{I \in \I_a}\left(\frac{1}{\mu_a(I)}\int_{I} |k_{a, \lambda}(x)|\,d\mu_a(x)\right)^{p} \\
&= \|k_{a,\lambda}\|^{2p}\mu_a(\{\lambda\})^{p}\sum_{I \in \I_a}\frac{\chi_{I}(\lambda)}{\mu_a(I)^{p}} \\
&\le c_p \|k_{\frac{a}{2},\lambda}\|^{2p}.
\end{align*}
Now let $\lambda$ be an arbitrary point in  $\C\setminus\supp\mu_a$. Then 
$k_{a,\lambda}(x)=- \frac{1}{2\pi i}\frac{1 - e^{-ia\bar\lambda}}{x -\bar\lambda}$ for all $x \in \supp\mu_a$. Thus, we need to estimate an oscillation of the function $x \mapsto \frac{1}{x -\bar\lambda}$ on the lattice $\Z_a$. Divide collection $\I_a$ from Section \ref{s1} into two parts: 
$$
\I_{a,1} = \{I \in \I_a\colon I = I_{a,j,k}, \;\Re\lambda \notin I_{a,j, k-1} \cup I_{a,j, k} \cup I_{a,j, k+1}\}, \;\; \I_{a,2} = \I_{a}\setminus\I_{a,1}.
$$
For an interval $I\in \I_{a,1}$ with center $x_c$, define the polynomial $P_I$ of degree $[\frac{1}{p}]$ by 
\begin{equation}\label{eq11}
\frac{1}{x-\bar\lambda} - P_I(x) = \frac{(x-x_{c})^{[\frac{1}{p}] + 1}}{(x-\bar\lambda)(\bar\lambda-x_{c})^{[\frac{1}{p}] + 1}}.
\end{equation}
Using this polynomial, we can estimate
\begin{equation}\label{eq4}
\osc\left(\frac{1}{x-\bar\lambda}, I, \mu_a, \bigl[\tfrac{1}{p}\bigr]\right) \le \sup_{x \in I}\left|\frac{(x-x_{c})^{[\frac{1}{p}] + 1}}{(x-\bar\lambda)(\bar\lambda-x_{c})^{[\frac{1}{p}] + 1}}\right| \le \frac{|I|^{[\frac{1}{p}] + 1}}{\dist(\lambda,I)^{[\frac{1}{p}] + 2}},
\end{equation}
where $|I|$ denotes the length of $I$. Since $I \in \I_{a,1}$, we have $\dist(\lambda, I) \ge |I|$, hence 
\begin{equation}\label{eq20}
\sum_{I \in \I_{a,1}} \osc\left(\frac{1}{\bar\lambda-x}, I, \mu_a, \bigl[\tfrac{1}{p}\bigr]\right)^{p} \le \sum_{I \in \I_{a,1}}\frac{1}{|I|^{p}} 
\le c_p \cdot a^p.
\end{equation} 
We also will need a more accurate estimate for the left hand side of the inequality above in the case where $|\Im\lambda|$ is large. For every $j\ge 0$, let $\I_{a,1}^{j}$ be the set of intervals $I_{a,j,k}$, $k \in \Z$, belonging to the family $\I_{a,1}$. We have 
\begin{align*}
\sum_{I \in \I_{a,1}^{j}} \left(\frac{|I|^{[\frac{1}{p}] + 1}}{\dist(\lambda,I)^{[\frac{1}{p}] + 2}}\right)^p = 
&\sum_{I \in \I_{a,1}^{j}} \left(\frac{|I|^{[\frac{1}{p}] + 1}}{\bigl(|\Im\lambda|^2+\dist(\Re\lambda,I)^{2}\bigr)^{([\frac{1}{p}] + 2)/2}}\right)^p\\
&\le c_p \left(\frac{a}{2^{j}}\right)^{p} \sum_{m\ge 1} \left(\frac{1}{(\frac{a}{2^j})^2|\Im\lambda|^2 + m^2}\right)^{\frac{1}{2}[\frac{1}{p}]p + p} \\
&\le c_p\left(\frac{a}{2^{j}}\right)^{p}\gamma_j^{1-[\frac{1}{p}]p - 2p},
\end{align*}
where $\gamma_j = \max\bigl(1,\frac{a}{2^j}|\Im\lambda|\bigr)$. Indeed, the last inequality follows from elementary estimates 
$$
\sum_{m=1}^{\infty} m^{-1 - 2p} < \infty, \qquad \int_{1}^{\infty}\frac{dx}{(r^2 + x^2)^{s}} \le c_s r^{1-2s},
$$
where $r>0$, and the constant $c_s$ depends on $s>1/2$. Put $N_\lambda = \bigl[\log_2(a|\Im\lambda|)\bigr]$ if $a|\Im\lambda| \ge 2$ and $N_\lambda = 0$ otherwise. Note that $\tilde p = -1+[\frac{1}{p}]p+p$ is a positive number. It follows 
\begin{align*}
\sum_{I \in \I_{a,1}} \osc\left(\frac{1}{\bar\lambda-x}, I, \mu_a, \bigl[\tfrac{1}{p}\bigr]\right)^{p} 
\le &c_p\sum_{j=0}^{\infty}\left(\frac{a}{2^{j}}\right)^{p}\gamma_j^{1-[\frac{1}{p}]p - 2p}\\
\le &c_p a^{-\tilde p}|\Im\lambda|^{-\tilde p-p}\sum_{j=0}^{N_\lambda}2^{\tilde p j} + 
c_p \sum_{j=N_\lambda}^{\infty}\frac{a^p}{2^{pj}}\\ 
\le &\frac{c_p}{|\Im\lambda|^p}.
\end{align*}
Combining the last estimate with \eqref{eq20}, we get
$$
\sum_{I \in \I_{a,1}} \osc\left(\frac{1}{\bar\lambda-x}, I, \mu_a, \bigl[\tfrac{1}{p}\bigr]\right)^{p} \le c_p \min\left(a^p, \frac{1}{|\Im\lambda|^{p}}\right).
$$
Now consider the family $\I_{a,2} = \I_{a,21} \cup \I_{a,22}$,  
$$
\I_{a,21} = \{I \in \I_{a,2}\colon |I| \le |\Im\lambda|\}, \quad \I_{a,22} = \{I \in \I_{a,2}\colon |I| > |\Im\lambda|\}.
$$
For an interval $I \in \I_{a,21}$ we use the polynomial $P_I$ defined by \eqref{eq11}. Then formula~\eqref{eq4} implies 
$$
\sum_{I \in \I_{a,21}} \osc\left(\frac{1}{\bar\lambda-x}, I, \mu_a, \bigl[\tfrac{1}{p}\bigr]\right)^{p} \le 
\sum_{I \in \I_{a,21}} \left(\frac{|I|^{[\frac{1}{p}]+1}}{|\Im\lambda|^{[\frac{1}{p}]+2}}\right)^{p} \le \frac{c_p}{|\Im\lambda|^p}.
$$
Note that if $|\Im\lambda| < \frac{2\pi}{a}$, the set $\I_{a,21}$ is empty. This shows that we can write 
$$
\sum_{I \in \I_{a,21}} \osc\left(\frac{1}{\bar\lambda-x}, I, \mu_a, \bigl[\tfrac{1}{p}\bigr]\right)^{p} \le c_p \min\left(a^p, \frac{1}{|\Im\lambda|^p}\right).
$$ 
For $I \in \I_{a,22}$ we put $P_I = 0$. Denote by $x_0$ the nearest point to $\lambda$ in $\supp\mu_a$, and set $I' = I\setminus\{x\in \R: |x - \Re\lambda|<\pi/a\}$. We have
\begin{align*}
\frac{1}{\mu_a(I)}\int_{I}\left|\frac{1}{x-\bar\lambda}\right|\,d\mu_a(x) 
&\le \frac{\mu_a(\{x_0\})}{\mu_a(I)|x_0-\bar\lambda|} + \frac{1}{\mu_a(I)}\int_{I'}\frac{dx}{|x-\bar\lambda|} \\
&\le \frac{c}{a|I||x_0-\bar\lambda|} + \frac{c}{|I|}\int_{\pi a^{-1}}^{|I|}\frac{dx}{\sqrt{x^2 + |\Im\lambda|^2}}\\
&\le \frac{c}{a|I||x_0-\bar\lambda|} + \frac{c}{|I|}\min\!\left(\log\frac{a|I|}{\pi}, \log\left(\frac{|I|}{|\Im\lambda|}+1\right)\!\right).
\end{align*}
Using estimates
$$
\sum_{I \in \I_{a,2}} \frac{1}{|I|^p} \le c_p a^p, \;\; \sum_{I \in \I_{a,2}} \left(\frac{\log a |I|}{|I|}\right)^p \le c_pa^p, \;
\;\sum_{I \in \I_{a,22}} \left(\frac{1}{|I|}\log\frac{|I|}{|\Im\lambda|}\right)^p \le \frac{c_p}{|\Im\lambda|^p},
$$
we see that
$$
\sum_{I \in \I_{a,22}}\osc\left(\frac{c_p}{\bar\lambda-x}, I, \mu_a, \bigl[\tfrac{1}{p}\bigr]\right)^{p} \le  \frac{c_p}{|x_0-\bar\lambda|^p} + c_p \min\left(a^p,\; \frac{1}{|\Im\lambda|^{p}}\right).
$$
Eventually, we obtain
$$
\left\|\frac{1}{x -\bar\lambda}\right\|_{\Bb_p(a,\osc)}^{p} \le \frac{c_p}{|x_0-\bar\lambda|^p} + c_p\min\left(a^p, \; \frac{1}{|\Im\lambda|^{p}}\right). 
$$
It follows that 
$$
\|k_{a,\lambda}\|_{\Bb_p(a,\osc)}^{p} \le c_p(1 + e^{-a\Im\lambda})^p\min\left(a^p, \; \frac{1}{|\Im\lambda|^{p}}\right) 
+ c_p \left|\frac{1 - e^{-ia\bar\lambda}}{x_0 - \lambda}\right|^p 
\le c_p \|k_{\frac{a}{2},\lambda}\|^{2p}, 
$$ 
which is the desired estimate. \qed

\bigskip

Let $\Cc_0(\R)$ denote the set of all continuous functions on $\R$ tending to zero at infinity. For completeness, we include the proof of the following known lemma.

\begin{Lem}\label{l7}
Let $0<p\le 1$, $a>0$. For every function $f \in \Bb_p(\osc,a)$ there exists a function $F \in \Bb_p(\R)$ such that $F = f$ on $\Z_a$, and
$$
\|F\|_{\Bb_p(\R)} \le c_p \|f\|_{\Bb_p(\osc,a)},
$$
where the constant $c_p$ depends only $p$.
\end{Lem}
\beginpf For $k \in \Z$ put $I_k = \left[\frac{2\pi}{a}[\frac{1}{p}]k, \frac{2\pi}{a}[\frac{1}{p}](k+1)\right]$. Interiors of intervals $I_k$ are disjoint and every set $I_k\cap\Z_a$ contains $[\frac{1}{p}] + 1$ points. On every $I_k$ define the polynomial $P_k$ of degree at most $[\frac{1}{p}]$ such that $P_k(x) = f(x)$ for all $x \in I_k \cap \Z_a$. Next, set $F(x) = P_k(x)$ for $x \in I_{k}$. 
We claim that the function $F$ is in~$\Bb_p(\R)$. To check this, let us take an interval 
$J_{j,k} = \left[\frac{2\pi}{a}[\tfrac{1}{p}]  k \cdot2^{j}, \frac{2\pi}{a}[\tfrac{1}{p}] (k+1) \cdot 2^{j}\right]$
with $k, j \in \Z$. In the case where $j < 0$ we clearly have $\osc(F,J_{j,k},m,[\tfrac{1}{p}]) = 0$ because the function $F$ is a polynomial of degree at most $[\frac{1}{p}]$ on $I$. Hence, we can assume that $J = J_{j,k}=I_{\ell} \cup \ldots \cup I_{\ell+N}$ for some $\ell \in \Z$ and $N \ge 1$. Consider the polynomial $P_{J}$ of degree at most $[\frac{1}{p}]$ such that 
$$
\osc\left(f,J,\mu_a,[\tfrac{1}{p}]\right) = \frac{1}{\mu_{a}(J)}\int_{J}|f(x) - P_J(x)|\,d\mu_a(x).
$$     
We have
\begin{multline*}
\frac{1}{|J|}\int_{J}|F(x) - P_J(x)|\,dx = \frac{1}{|J|}\sum_{s=0}^{N} \int_{I_{\ell+s}}|P_{\ell + s}(x) - P_J(x)|\,dx \le \\
\le \frac{c_p}{|J|}\sum_{s=0}^{N} \int_{I_{\ell+s}}|P_{\ell + s}(x) - P_J(x)|\,d\mu_a(x) \le c_p  \osc\left(f,I,\mu_a,[\tfrac{1}{p}]\right),
\end{multline*}
where we used the fact that 
$$
\int_{I_\ell} |P(x)|\,dx  \le c_p \int_{I_\ell} |P(x)|\,d\mu_a(x)
$$
for every interval $I_\ell$, $\ell \in \Z$, and every polynomial $P$ of degree at most $[\frac{1}{p}]$. It follows that 
$$
\|F\|_{\Bb_p(\R,m,\osc)}^{p} \le c_p^p\sum_{j,k}\osc\left(f,J_{j,k},\mu_a,[\tfrac{1}{p}]\right)^p \le c_p^p \|f\|_{\Bb_p(\osc, a)}^{p},
$$
and hence $F$ belongs to the space $\Bb_{p,p}^{1/p}(\R,dx,\osc) = \Bb_{p}(\R)$, as required. \qed

\bigskip

\noindent{\bf Proof of Proposition \ref{P4}.} Consider a function $F$ of the form 
$$
F = \sum_{\lambda \in \Lambda} c_\lambda \frac{k_{a,\lambda}}{\|k_{\frac{a}{2},\lambda}\|^2}, \qquad 
\sum_{\lambda \in \Lambda}|c_{\lambda}|^{p} < \infty.
$$ 
Since $0<p\le 1$ and $|k_{a,\lambda}(x)| \le c\|k_{\frac{a}{2}, \lambda}\|^{2}$ for every $\lambda \in \C$ and $x \in \R$, the series above converges absolutely to a function from $\Cc_0(\R)$ by the Lebesgue dominated convergence theorem. By Lemma~\ref{l17}, the restriction of $F$ to $\Z_a$ (to be denoted by $f$) is in $\Bb_p(a,\osc)$ and $\|f\|_{\Bb_p(a,\osc)}^{p} \le c_p \sum_{\lambda \in \Lambda}|c_\lambda|^{p}$ for a constant $c_p$ depending only on $p$. 

\medskip

Conversely, take $f \in \Bb_p(a,\osc)$ and find a function $\tilde F \in \Bb_p(\R)$ such that $\tilde F = f$ on $\Z_a$, see Lemma \ref{l7}. Applying Theorem 2.10 from \cite{Roch85} to analytic and anti-analytic parts of $\tilde F$, we obtain the representation
$$
\tilde F = q-\frac{1}{2\pi i}\sum_{\lambda \in \U_{\eps}} \tilde c_{\lambda} \frac{\Im\lambda}{x-\bar\lambda}, \qquad x \in \R,
$$
where the coefficients $\tilde c_k \in \C$ are such that $\sum|\tilde c_{\lambda}|^p \le c_p\|\tilde F\|_{\Bb_p(\R)}^{p}$, and $q$ is a polynomial of degree at most $[\frac{1}{p}]$. Now consider the function 
$$
F = \sum_{\lambda \in \U_{\eps}} c_{\lambda} \frac{k_{\lambda,a}}{\|k_{\frac{a}{2},\lambda}\|^2}, \qquad 
c_{\lambda} = \tilde c_{\lambda} \frac{\Im\lambda\cdot \|k_{\frac{a}{2},\lambda}\|^2}{1-e^{-ia\bar\lambda}}. 
$$
Observe that $|c_\lambda| \le |\tilde c_\lambda|$ for all $\lambda \in \U_{\eps}$ and $f = q+F$ on $\Z_a$. We need to replace the set $\U_\eps$ above to the set $\Lambda_{\eta a,\eps}$ from Proposition \ref{p1}. Since $k_{\frac{a}{2},\lambda} = e^{\frac{iaz}{4}}e^{-\frac{ia\bar\lambda}{4}}\rho_{\frac{a}{4},\lambda}$, we have $\|k_{\frac{a}{2},\lambda}\|^{2} = e^{-\frac{a\Im\lambda}{2}}\|\rho_{\frac{a}{4},\lambda}\|^{2}$ and
$$
e^{-\frac{iax}{2}}F = \sum_{\lambda \in \U_{\eps}} c_{\lambda} e^{-\frac{ia\bar\lambda}{2}} \frac{\rho_{{a/2},\lambda}}{\|k_{a,\lambda}\|^2} 
= \sum_{\lambda \in \U_{\eps}} c_{\lambda} e^{-\frac{ia\Re\lambda}{2}}\frac{\rho_{a/2,\lambda}}{\|\rho_{a/4,\lambda}\|^2}.
$$
From the beginning of the proof of Proposition \ref{p1} we see that the Toeplitz operator on $\pw_{a/4}$ with symbol $e^{-\frac{iax}{2}}F$ belongs to the class $\Ss^p(\pw_{a/4})$. It follows that  
$$
e^{-\frac{iax}{2}}F = \sum_{\lambda \in \Lambda_{\eta a, \eps}} d_{\lambda} \frac{\rho_{a/2,\lambda}}{\|\rho_{a/4,\lambda}\|^2}, \qquad
\sum_{\lambda \in \Lambda_{\eta a,\eps}}|d_\lambda|^p \le c_p \sum_{\lambda \in \U_{\eps}}|c_\lambda|^p.
$$ 
This yields the required representation for $F$, 
$$
F = \sum_{\lambda \in \Lambda_{\eta a,\eps}} c_{\lambda} \frac{k_{a,\lambda}}{\|k_{\frac{a}{2},\lambda}\|^2},
\qquad \sum_{\lambda \in \Lambda_{\eta a,\eps}}|c_\lambda|^p \le c_p \|f\|_{\Bb_p(a,\osc)},
$$
with some new coefficients $c_\lambda$. Since $\sum_{\lambda} |c_\lambda| < \infty$, the function $G = e^{\frac{-iaz}{2}}F$ is an entire function of exponential type at most $a/2$ such that $\lim_{x \to \pm\infty} |G(x)| = 0$. In particular, it is uniquely determined by values on $\Z_a$. This proves uniqueness in Proposition \ref{P4}. \qed 

\bigskip

\section{Proof of Theorem \ref{t1} for $0<p\le 1$}\label{s5}

\noindent{\bf Proof of Theorem \ref{t1} ($0<p\le 1$).} Let $\phi \in L^\infty(\R)$ be a function on $\R$ such that the operator $T_{\phi}$ is in $\Ss^p(\pw_a)$, and let $\phi_{st} = \F^{-1}\chi_{(-2a,2a)}\F\phi$ be the standard symbol of $T_{\phi}$. By Proposition \ref{p1} and Proposition~\ref{P4}, we have $e^{2iax}\phi_{st} \in \Bb_p(4a, \osc)$ and moreover, 
$\|e^{2iax}\phi_{st}\|_{\Bb_p(4a,\osc)} \le c_p \|T_{\phi}\|_{\Ss^p}$ for a constant $c_p$ depending only on $p$. 

\medskip

Conversely, assume that the restriction of the function $e^{2iax}\phi_{st}$ to $\Z_{4a}$ belongs to the space $\Bb_p(4a, \osc)$. By Proposition \ref{P4}, there exists a function $F$ and a polynomial $q$ of degree at most $[\frac{1}{p}]$ such that $q + F= e^{2iax}\phi_{st}$ on $\Z_{4a}$ and 
\begin{equation}\label{eq17}
F = \sum_{\lambda \in \Lambda_{\eta a,\eps}} c_\lambda \frac{k_{4a,\lambda}}{\|k_{2a,\lambda}\|^2}=
e^{2iax}\sum_{\lambda \in \Lambda_{\eta a,\eps}} c_\lambda e^{-2ia\Re\lambda}\frac{\rho_{2a,\lambda}}{\|\rho_{a,\lambda}\|^2}
\end{equation}
for some $c_\lambda \in \C$ such that $\sum|c_\lambda|^p \le c_p \|e^{2iax}\phi_{st}\|^{p}_{\Bb_{p}(4a,\osc)}$. We claim that $T_{\tilde\phi} = T_{\phi}$ on $\pw_a$, where $\tilde\phi = e^{-2iax}(q+F)$. Indeed, the entire function $z \mapsto \tilde\phi - \phi_{st}$ has exponential type at most $2a$, vanishes on $\Z_{4a}$, and satisfies a polynomial estimate on $\R$. Hence $\tilde\phi-\phi_{st} = \tilde q\sin (2az)$ for all $z \in \C$ and a polynomial $\tilde q$. Thus, we have $T_\phi = T_{\phi_{st}} = T_{\tilde \phi}$. It remains to use formula~\eqref{eq17} and  Proposition \ref{p1}. The theorem is proved. \qed

\bigskip

\section{Discrete Hilbert transform commutators. Proof of Theorem \ref{t2}}\label{s4} 
Recall that $\mu_a = \frac{2\pi}{a}\sum_{x \in \Z_a} \delta_{x}$ is the scalar multiple of the counting measure on the lattice  $\Z_{a} = \bigl\{\frac{2\pi}{a}k, \; k \in \Z\bigr\}$. The discrete Hilbert transform $H_{\mu_a}$ on $L^2(\mu_a)$ is defined by 
$$
H_{\mu_a}: f \mapsto \frac{1}{\pi}\dashint_{\Z_a}\frac{f(t)}{x-t}\,d\mu_a(t),
$$ 
and its commutator $C_\psi = M_\psi H_{\mu_a} - H_{\mu_a}M_\psi$ with the multiplication operator $M_\psi: f \mapsto \psi f$ on $L^2(\mu_a)$ by
$$
C_{\psi}: f \mapsto  \frac{1}{\pi} \dashint_{\Z_a} \frac{\psi(x) - \psi(t)}{x -t} f(t)\,d\mu_a(t), \qquad x \in \supp\mu_{a}.
$$
It is well-known that the operator $H_{\mu_a}$ admits the bounded extension from the dense subset $\mathcal{G}$ of $L^2(\mu_a)$ of finitely supported bounded functions to the whole space $L^2(\mu_a)$.  A possible way to define the operator $C_\psi$ on $L^2(\mu_a)$ for any symbol $\psi$ on $\Z_a$ is to consider its bilinear form on elements from the dense subset $\mathcal{G}\times\mathcal{G}$ of~$L^2(\mu_a) \times L^2(\mu_a)$. We will also deal with the operators $\tilde C_{\psi}: L^2(\mu_{\frac{a}{2}}) \to L^2(\nu_{\frac{a}{2}})$ defined by
$$
\tilde C_{\psi}: f \mapsto  \frac{1}{\pi} \int_{\Z_a} \frac{\psi(x) - \psi(t)}{x -t} f(t)\,d\mu_{\frac{a}{2}}(t), \qquad x \in \supp\nu_{\frac{a}{2}},
$$
where the measure $\nu_{\frac{a}{2}} = \frac{4\pi}{a}\sum_{x \in \Z_{\frac{a}{2}}}\delta_{x + \frac{2\pi}{a}}$ is supported on the lattice 
$\frac{2\pi}{a} + \Z_{\frac{a}{2}}$. It can be shown that for $1\le p\le\infty$ the operator $C_{\psi}:L^2(\mu_a) \to L^2(\mu_a)$ is in $\Ss^p$ if and only if the operator $\tilde C_{\psi}: L^2(\mu_{\frac{a}{2}}) \to L^2(\nu_{\frac{a}{2}})$ is in $\Ss^p$. As we will see, for $0<p<1$ we may have $C_{\psi} \notin \Ss^p(L^2(\mu_a))$ for a function $\psi$ on $\Z_a$ such that the operator $\tilde C_{\psi}: L^2(\mu_{\frac{a}{2}}) \to L^2(\nu_{\frac{a}{2}})$ is in $\Ss^p$.  

\medskip 

The discrete Hilbert transform commutators were investigated in details in paper~\cite{Roch87}. In particular, it was proved in \cite{Roch87} that $C_{\psi}$ is bounded on $L^2(\mu_a)$ if and only if its symbol $\psi$ belongs to the discrete $\bmo(\Z_a)$ space of functions $f$ on $\Z_a$ such that $\sup_{I\in \I_a}\osc(f,I,\mu_a,0) <\infty$, where $\I_a = \{I_{a,j,k},\; j,k\in \Z, \; j\ge0\}$ is the collection of intervals defined in Section \ref{s1}. Another result from \cite{Roch87} says that $C_{\psi}$ is compact on $L^2(\mu_a)$ if and only if $\psi \in \cmo(\Z_a)$, that is, $\lim_{k \to \pm\infty}\osc(\psi,I_{a,j,k},\mu_a,0) = 0$ for every $j\ge 0$ and $\lim_{j\to+\infty}\osc(\psi,J_j,\mu_a,0) = 0$  for any sequence of intervals $J_{j} \subset\R$ of length $j$ with common center. Finally, the operator $C_\psi$ belongs to $\Ss^p(L^2(\Z_a))$ for $1<p<\infty$ if and only if $\psi \in \Bb_p(a,\osc)$, moreover, we have $C_{\psi} \in \Ss^1(L^2(\mu_a))$ for every  $\psi \in \Bb_1(a,\osc)$. See Theorem 6.2 in \cite{Roch87} and Theorem 4 in \cite{Tor98} for the proof of these results. It was an open question stated in Section 7 of \cite{Roch87} whether $C_{\psi} \in \Ss^p(L^2(\mu_a))$ is equivalent to $\psi \in \Bb_p(a,\osc)$ for all positive $p$ (in particular, for $p=1$). 
Theorem~\ref{t2} gives the affirmative answer to this question for $p=1$. On the other hand, for $0<p<1$ we show that there exists symbols $\psi \in \Bb_p(a,\osc)$ such that $C_{\psi} \notin \Ss^p(L^2(\mu_a))$. In fact, the following modification of Theorem \ref{t2} holds true. 
\addtocounter{Thm}{-1}
\renewcommand{\theThm}{\arabic{Thm}$'$} 
\begin{Thm}\label{t2p}
Let $0<p\le 1$. The operator $\tilde C_{\psi}: L^2(\mu_{\frac{a}{2}}) \to L^2(\nu_{\frac{a}{2}})$ belongs to the class $\Ss^p$ if and only if 
 $\psi \in \Bb_p(a,\osc) \cap L^\infty(\Z_{a})$. Moreover, the quasi-norms 
$\|\tilde C_{\psi}\|_{\Ss^p}$ and $\|\psi\|_{\Bb_p(a,\osc)}$ are comparable with constants depending only on $p$.
\end{Thm} 
\renewcommand{\theThm}{\arabic{Thm}} 
For the proof we need a result on unitary equivalence of discrete Hilbert transform commutators to some truncated Hankel operators.   
Given a positive number $a>0$, we denote by $\pw_{[-a,0]}$ the Paley-Wiener space of functions in $L^2(\R)$ with Fourier spectrum in the interval $[-a,0]$. Define the truncated Hankel operator $\Gamma_{\psi}: \pw_{[0,a]} \to \pw_{[-a,0]}$ with symbol $\psi \in L^\infty(\R)$ by
$$
\Gamma_{\psi}: f \mapsto P_{[-a,0]}(\psi f), \qquad f \in \pw_{[0,a]},
$$ 
where $P_{[-a,0]}$ stands for the projection in $L^2(\R)$ to the subspace $\pw_{[-a,0]}$. It is easy to see that $\Gamma_{\psi}$ is completely determined by its standard symbol $\psi_{st,2a} = \F^{-1}\chi_{(-2a,0)}\F\psi$, that is, $\Gamma_{\psi}f = \Gamma_{\psi_{st,a}}f$ for all functions $f\in \pw_{[0,a]}$ such that $\sup_{x\in \R}|xf(x)|<\infty$. Clearly, such functions form a dense subset in $\pw_{[0,a]}$.  

\medskip

It is known that the embedding operator $V_{\mu_a}: \pw_{[0,a]} \to L^2(\mu_a)$ taking a function $f \in \pw_{[0,a]}$ into its restriction to $\Z_a$ is unitary. The same is true for the embedding operator $\tilde V_{\nu_a}: \pw_{[-a,0]} \to L^2(\nu_a)$. A general version of the following result is Lemma 4.2 of \cite{Bes15}.

\medskip

\begin{Lem}\label{l12}
Let $a>0$, $0<p\le 1$, and let $\psi \in L^{\infty}(\Z_{2a})$. 
Then there exists an entire function $\Psi$ such that $\Psi = \psi$ on $\Z_{2a}$, 
$|F(x)| \le c\log(e+|x|)$ for all $x \in \R$, and the Fourier spectrum of $F$ is contained in the interval $[-2a,0]$. Moreover, we have
\begin{equation}\label{eq22}
\tilde V_{\nu_a}\Gamma_{\Psi}V_{\mu_a}^{-1} = -i \tilde  C_{\psi}.
\end{equation}
for the operators $\Gamma_{\Psi}: \pw_{[0,a]}\to\pw_{[-a,0]}$ and $\tilde C_{\psi}: L^2(\mu_{a}) \to L^2(\nu_{a})$.
\end{Lem}
\beginpf Existence of such a function $\Psi$ follows from a general theory of entire functions, see, e.g., Theorem 1 in Section 21.1 of \cite{Levin} and Problem~1 after its proof. In order to prove formula \eqref{eq22}, take a pair of functions $f \in L^2(\mu_{a})$, $g \in L^2(\nu_{a})$ with finite support. Consider the functions $F,G$ in $\pw_{[0,a]}$ such that $F = V_{\mu_a}^{-1}f$, $\bar G = \tilde V_{\nu_a}^{-1}g$. It is easy to see that $\int_{\R}|\Psi FG|\,dx <\infty$ and hence the bilinear form of $\Gamma_{\Psi}$ is correctly defined on functions $F$, $\bar G$. We have
\begin{multline*}
(\tilde V_{\nu_a}\Gamma_{\Psi}V_{\mu_a}^{-1}f,g)_{L^2(\R)} 
= (\Gamma_{\Psi} F, \bar G)_{L^2(\R)} = (FG, \bar\Psi)_{L^2(\R)} =\\  
= (V_{\mu_{2a}}FG, V_{\mu_{2a}}\bar \Psi)_{L^2(\mu_{2a})} = 
\tfrac{1}{2}(Fg, \bar\psi)_{L^2(\nu_{a})} + \tfrac{1}{2}(fG,\bar\psi)_{L^2(\mu_{a})}.
\end{multline*} 
For every point $x \in \frac{\pi}{a} + \Z_{a}$ we have 
$$
F(x) = (V_{\mu_a} F, V_{\mu_a} k_{x,a})_{L^2(\mu_a)} = \frac{2}{\pi i}\int_{\R}\frac{f(t)}{t-x}\,d\mu_a(t), \qquad x \in \frac{\pi}{a} + \Z_{a}.
$$
Analogously, $G(t) = \frac{2}{\pi i}\int_{\R}\frac{\ov{g(x)}}{x-t}\,d\nu_a(x)$ for all $t \in \Z_{a}$. Using these formulas, we get
\begin{align*}
(\tilde V_{\nu_a}\Gamma_{\Psi}V_{\mu_a}^{-1}f,g)_{L^2(\R)} 
&= \frac{1}{\pi i}\int_{\R}\frac{\psi(x) - \psi(t)}{x-t}f(t)\ov{g(x)}\,d\mu_a(t)\,d\nu_a(x) \\ 
&= -i (\tilde C_{\psi} f, g)_{L^2(\nu_{a})}.
\end{align*}
The lemma follows. \qed

\medskip

\noindent {\bf Proof of Theorem \ref{t2p}.} Let $\psi$ be a function on the lattice $\Z_a$ such that the operator $\tilde C_{\psi}: L^2(\mu_{\frac{a}{2}}) \to L^2(\nu_{\frac{a}{2}})$ belongs to the class $\Ss^p$. Consider the sequence of points $x_k = \frac{2\pi}{a}k$, $k\in \Z$. Since $0<p\le 1$, we have
$$
\sum_{k \in \Z} |\psi(x_{2k}) - \psi(x_{2k+1})| = \frac{a}{8}\sum_{k \in \Z} |(\tilde C_{\psi} \delta_{x_{2k}}, \delta_{x_{2k+1}})_{L^2(\nu_{\frac{a}{2}})}| <\infty.
$$ 
Hence, the function $\psi$ is bounded on $\Z_{a}$. Using Lemma \ref{l12}, we can find an entire function $\Psi$ such that $\Psi = \psi$ on $\Z_{a}$, 
$|\Psi(x)| \le c\log(e+|x|)$ for all $x \in \R$, the Fourier spectrum of $\Psi$ is contained in $[-a,0]$, and relation \eqref{eq22} holds for the operators $\Gamma_{\Psi}: \pw_{[0,\frac{a}{2}]}\to\pw_{[-\frac{a}{2},0]}$ and $\tilde C_{\psi}: L^2(\mu_{\frac{a}{2}}) \to L^2(\nu_{\frac{a}{2}})$.
In particular, we have $\Gamma_{\Psi} \in \Ss^p$. Denote by $M$ the multiplication operator on $L^2(\R)$ by the function $e^{\frac{iax}{2}}$. Let $T_{e^{\frac{iax}{2}}\Psi}$ be the Toeplitz operator on $\pw_{\frac{a}{4}}$ with standard symbol $e^{\frac{iax}{2}}\Psi$. Observe that 
\begin{equation}\label{eq21}
T_{e^{\frac{iax}{2}}\Psi}f = M\Gamma_{\Psi}Mf,
\end{equation}
for every function function $f \in \pw_{\frac{a}{4}}$ such that $\sup_{x\in \R} |xf(x)| < \infty$. Since $M$ maps unitarily $\pw_{\frac{a}{4}}$ onto $\pw_{[0,\frac{a}{2}]}$ and $\pw_{[-\frac{a}{2},0]}$ onto $\pw_{\frac{a}{4}}$, the operator  $T_{e^{\frac{iax}{2}}\Psi}$ belongs to $\Ss^p(\pw_{\frac{a}{4}})$.  In particular, there exists a function $\phi \in L^\infty(\R)$ such that $T_{\phi} = T_{e^{\frac{iax}{2}}\Psi}$ and
$\phi_{st} = e^{\frac{iax}{2}}\Psi + c_1e^{-i \frac{a}{2} x} + c_2e^{i \frac{a}{2} x}$ for some constants $c_1$, $c_2$. Since $e^{\frac{iax}{2}}\phi_{st}$ coincides with $\psi + c_1 + c_2$ on $\Z_a$, we have $\psi \in \Bb_p(a,\osc)$ by Theorem~\ref{t1}. Moreover, the quasi-norm $\|\tilde C_{\psi}\|_{\Ss^p}$ is comparable to $\|\psi\|_{\Bb_p(a,\osc)}$ with constants depending only on $p \in (0,1]$. 

\medskip

Conversely, suppose that $\psi \in \Bb_p(a,\osc)\cap L^\infty(\Z_a)$. Using Lemma \ref{l12} again, we find an entire function $\Psi$ such that $\Psi = \psi$ on $\Z_{a}$, $|\Psi(x)| \le c\log(e+|x|)$ for all $x \in \R$, the Fourier spectrum of $\Psi$ is contained in $[-a,0]$, and relation \eqref{eq22} holds for the operators $\Gamma_{\Psi}: \pw_{[0,\frac{a}{2}]}\to\pw_{[-\frac{a}{2},0]}$ and $\tilde C_{\psi}: L^2(\mu_{\frac{a}{2}}) \to L^2(\nu_{\frac{a}{2}})$. Since $\psi \in L^\infty(\Z_a)$, the operators $\tilde C_{\psi}$ and $\Gamma_{\Psi}$ are bounded. Let $\Psi_{st,a}$ be the standard symbol of the operator $\Gamma_\Psi$. Note that $\Psi_{st,a}(x) = \Psi(x) + q(x)$ for all $x \in \Z_a$ and a polynomial $q$ of degree at most one. In particular, we have $\Psi_{st,a} \in \Bb_p(a, \osc)$. By Theorem \ref{t1}, the operator $T_{e^{\frac{iax}{2}}\Psi_{st,a}}$ on $\pw_{\frac{a}{4}}$ is in $\Ss^p$, hence $\Gamma_{\Psi} \in \Ss^p$ by formula~\eqref{eq21}. It follows that the operator $\tilde C_{\psi}$ is in $\Ss^p$ as well, and, moreover, we have the estimate 
$$
\|\tilde C_{\psi}\|_{\Ss^p} = \|\Gamma_{\Psi}\|_{\Ss^p} = \left\|T_{e^{\frac{iax}{2}}\Psi_{st,a}}\right\|_{\Ss^p} \le c_p \|\Psi_{st,a}\|_{\Bb_p(a,\osc)} = c_p \|\psi\|_{\Bb_p(a,\osc)},
$$ 
for a constant $c_p$ depending only on $p$. The theorem is proved. \qed

\bigskip

\noindent {\bf Proof of Theorem \ref{t2}.} Let $\psi$ be a function on the lattice $\Z_a$ such that we have $C_{\psi} \in \Ss^1(L^2(\mu_a))$. Then the operator $\tilde C_\psi: L^2(\mu_{\frac{a}{2}}) \to L^2(\nu_{\frac{a}{2}})$ is of trace class as well and $\|\psi\|_{\Bb_1(a,\osc)} \le c_1\|\tilde C_\psi\|_{\Ss^1(L^2(\mu_a))} \le c_1\|C_\psi\|_{\Ss^1(L^2(\mu_a))}$ by Theorem \ref{t2p}. 

\medskip

Conversely, suppose that $\psi \in \Bb_1(a,\osc)\cap L^\infty(\Z_a)$. By Lemma \ref{l7}, we can find a function $\Psi \in \Bb_1(\R) \cap L^\infty(\R)$ such that $\Psi = \psi$ on $\Z_a$ and $\|\Psi\|_{\Bb_1(\R)} \le c_1 \|\psi\|_{\Bb_1(\osc,a)}$. Denote $\psi_\lambda: t \mapsto \frac{|\Im\lambda|^2}{(t - \bar\lambda)^2}$ for $\lambda \in \C$. Let us apply Theorem~2.10 in \cite{Roch85} to analytic and anti-analytic parts of $\Psi$: find numbers $c$, $c_\lambda$ such that 
$\sum_{\lambda \in \U_{\eps}}|c_{\lambda}| \le c_1\|\Psi\|_{\Bb_1(\R)}$ and
$$
\psi(x) = \Psi(x) = c+ \sum_{\lambda \in \U_\eps} c_{\lambda} \psi_{\lambda}(x), \qquad x \in \Z_a. 
$$ 
We claim that for every $\lambda \in \U_\eps$ the commutator $C_{\psi_\lambda}$ belongs to the trace class and $\|C_{\psi_\lambda}\|_{\Ss^1} \le c_1 (1+ a)$ for a constant $c_1$ do not depending on $\lambda$. Clearly, this will yield the desired estimate $\|C_{\psi}\|_{\Ss^1} \le c_1 (1+a)\|\psi\|_{\Bb_1(a,\osc)}$. 
We have
\begin{align*}
\frac{\psi_\lambda(x) - \psi_\lambda(t)}{x -t} = -\frac{|\Im\lambda|^2}{(x - \bar\lambda)^2(t - \bar\lambda)} - \frac{|\Im\lambda|^2}{(x - \bar\lambda)(t - \bar\lambda)^2}. 
\end{align*}
Denote by $K_{\psi_\lambda}$ the integral operator on $L^2(\mu_a)$ with kernel $\frac{\psi_\lambda(x) - \psi_\lambda(t)}{x -t}$:
\begin{equation}\label{eq19}
(K_{\psi_\lambda}f)(x) = \int_{\Z_a}\frac{\psi_\lambda(x) - \psi_\lambda(t)}{x -t}f(t)\,dt = (C_{\psi_\lambda}f)(x) + \frac{2|\Im\lambda|^2}{(x-\bar\lambda)^{3}}f(x).  
\end{equation}
Observe that the operator $K_{\psi_\lambda}$ has rank $2$ and 
$$
\|K_{\psi_\lambda}\|_{\Ss^p} \le 2|\Im\lambda|^2\cdot \left\|\frac{1}{(x - \bar\lambda)^2}\right\|_{L^2(\mu_a)} \left\|\frac{1}{x - \bar\lambda}\right\|_{L^2(\mu_a)}. 
$$
In the case where $\dist(\lambda, \Z_a) \ge \frac{\pi}{2a}$, the last expression could be estimated from above by 
$$
c_1 \left(\int_{\R}\frac{|\Im\lambda|\,dt}{t^2 + |\Im\lambda|^2}\int_{\R}\frac{|\Im\lambda|^3\,dt}{(t^2 + |\Im\lambda|^2)^{2}}\right)^{\frac{1}{2}} = c_1\left(\int_{\R}\frac{dt}{t^2 +1}\int_{\R}\frac{dt}{(t^2 +1)^2} \right)^{\frac{1}{2}}.
$$
Moreover, the singular numbers of the multiplication operator $f \mapsto \frac{|\Im\lambda|^2}{(x-\bar\lambda)^{3}}f$ are precisely 
$\frac{|\Im\lambda|^2}{|x-\bar\lambda|^{3}}$, $x \in \Z_a$, hence its norm in $\Ss^1(L^2(\mu_a))$ does not exceed
$$
\sum_{x \in \Z_a} \frac{|\Im\lambda|^2}{|x-\bar\lambda|^{3}} \le \sum_{x \in \Z_a} \frac{|\Im\lambda|^2}{(x^2+|\Im\lambda|^2)^{\frac{3}{2}}} \le c_1 a
$$
for a universal constant $c_1$. This tells us that $\|C_{\psi_\lambda}\|_{\Ss^p} \le c_1(1+a)$ for all $\lambda \in \U_\eps$ such that $\dist(\lambda, \Z_a) \ge \frac{\pi}{2a}$. Now consider the case where $\dist(\lambda, \Z_a) \le \frac{\pi}{2a}$. Let $x_\lambda$ be the nearest point to $\lambda$ in the lattice $\Z_a$. The function $\psi_\lambda$ belongs to $L^1(\mu_a)$ and 
\begin{align*}
\sum_{x \in \Z_a}|\psi_\lambda(x)|
&\le |\psi_\lambda(x_\lambda)|+ 2|\Im\lambda|^2 \sum_{k=1}^{\infty}\frac{1}{(\frac{2\pi}{a}k - \frac{\pi}{2a})^2}, \\
&\le \left|\frac{\Im\lambda}{\lambda - x_\lambda}\right|^2 +  2\left(\frac{a|\Im\lambda|}{2\pi}\right)^2\sum_{k=1}^{\infty}\frac{1}{(k - \frac{1}{4})^2}
\le c_1,
\end{align*} 
where the right hand side does not depend on $\lambda$. It follows that the operator $M_{\psi_\lambda}$ lies in $\Ss^1(L^2(\mu_a))$ and $\|M_{\psi_\lambda}\|_{\Ss^1} \le c_1$.  We also have 
$$
\bigl\|C_{\psi_\lambda}\bigr\|_{\Ss^p} = \bigl\|H_{\mu_a}M_{\psi_\lambda} - M_{\psi_\lambda} H_{\mu_a}\bigr\|_{\Ss^1} \le c_1,
$$ 
for another constant $c_1$, because the discrete Hilbert transform $H_{\mu_a}$ is bounded on~$L^2(\mu_a)$. This completes the proof. \qed

\bigskip

Remark that the second part of the proof of Theorem \ref{t2} is almost literal repetition of the corresponding part of the proof of Theorem 6.2 in \cite{Roch87}. However, the original argument in \cite{Roch87} has a gap: it does not involve the estimate of the $\Ss^1$-norm of the multiplication operator $f \mapsto \frac{|\Im\lambda|^2}{(x-\bar\lambda)^{3}}f$ from formula \eqref{eq19}. This technical place turns out to be crucial in the case $0<p<1$. More precisely, we have the following result.
\begin{Prop}\label{p3}
Let $0<p<1$ and let $a>0$. There exists a function $\psi \in \Bb_{p}(\Z_a)$ such that $C_{\psi} \notin \Ss^p(L^2(\mu_a))$. 
\end{Prop}
\beginpf Suppose that $C_{\psi} \in \Ss^p(L^2(\mu_a))$ for every $\psi \in \Bb_{p}(a,\osc)$. Then it is easy to see from the closed graph theorem that there exists a constant $c_{p,a}$ such that $\|C_{\psi}\|_{\Ss^p} \le c_{p,a}\|\psi\|_{\Bb_p(a,\osc)}$ for all $\psi \in \Bb_{p}(a,\osc)$. Take $\lambda \in \C^+$ such that $\Im\lambda \ge \frac{2\pi}{a}$ and consider the function  $\psi_\lambda: t \mapsto \frac{\Im\lambda}{t - \bar\lambda}$. Analogously to \eqref{eq19}, we have $K_{\psi_{\lambda}} = C_{\psi_{\lambda}} + M_\lambda$, where $K_{\psi_{\lambda}}$ is the integral operator with kernel 
$$
\frac{\psi_\lambda(x) - \psi_\lambda(t)}{x -t} = - \frac{\Im\lambda}{(x-\bar\lambda)(t-\bar\lambda)},
$$ 
and $M_\lambda: f \mapsto \frac{\Im\lambda}{(x-\bar\lambda)^2} f$ is the multiplication operator on $L^2(\mu_a)$ by $\frac{\Im\lambda}{(x-\bar\lambda)^2}$. Observe that $K_{\psi_\lambda}$ is the rank-one operator whose norm does not exceed 
$$
\Im\lambda\cdot\left\|\frac{1}{x-\bar\lambda}\right\|_{L^2(\mu_a)}^{2} \le c_p \int_{\R} \frac{\Im\lambda\,dt}{t^2 + (\Im\lambda)^2} 
= c_p \int_{\R} \frac{dt}{t^2 + 1}.
$$
It follows from our assumption and Lemma \ref{l17} that $\|M_\lambda\|_{\Ss^p} \le c_{p,a}$ for all $\lambda\in \C$ with $\Im\lambda \ge  \frac{2\pi}{a}$ and a universal constant $c_{p}$. On the other hand, we have
$$
\|M_\lambda\|_{\Ss^p}^{p} = \sum_{x \in \Z_a} \frac{(\Im\lambda)^p}{|x-\bar\lambda|^{2p}} \ge 
ac_p\int_{\R}\frac{(\Im\lambda)^p\,dx}{(x^2 + (\Im\lambda)^{2})^{p}} = ac_p (\Im\lambda)^{1-p}\int_{\R}\frac{dt}{t^2+1}.
$$ 
Since the right hand side is unbounded in $\lambda$, we get the contradiction. \qed

\bibliographystyle{plain} 
\bibliography{bibfile}
\enddocument